\documentclass[11pt,a4paper]{amsart}
\usepackage{graphics,array,amsmath,amssymb,epsf,graphicx}
\usepackage{esvect}

\newtheorem{thm}{Theorem}
\newtheorem{lem}[thm]{Lemma}

\newtheorem{theorem}{Theorem}
\newtheorem{ctheorem}{Theorem}

\theoremstyle{definition}
\newtheorem*{definition}{Definition}

\begin{document}
\title{Linearly Free Graphs}
\author{Youngsik Huh}
\author{Jung Hoon Lee}
\address{Department of Mathematics,
College of Natural Sciences, Hanyang University, Seoul 04763,
Korea} \email{yshuh@hanyang.ac.kr}
\address{Department of Mathematics and Institute of Pure and Applied Mathematics,
Chonbuk National University,
Jeonju 54896,
Korea} \email{junghoon@jbnu.ac.kr}




\begin{abstract}
In this paper we are interested in an intrinsic property of graphs which is derived from their embeddings into the Euclidean 3-space $\mathbb{R}^3$.
An embedding of a graph into $\mathbb{R}^3$ is said to be {\em linear}, if it sends every edge to be a line segment. And we say that an embedding $f$ of a graph $G$ into $\mathbb{R}^3$ is {\em free}, if $\pi_1(\mathbb{R}^3-f(G))$ is a free group.
Lastly a simple connected graph is said to be {\em linearly free} if every its linear embedding is free.
In 1980s it was proved that every complete graph is linearly free, by Nicholson \cite{Nicholson}.

In this paper, we develop Nicholson's arguments into a general notion, and establish a sufficient condition for a linear embedding to be free.
As an application of the condition we give a partial answer for a question: how much can the complete graph $K_n$ be enlarged so that the linear freeness is preserved and the clique number does not increase? And an example supporting our answer is provided.

As the second application it is shown that a simple connected graph of minimal valency at least $3$ is linearly free, if it has less than 8 vertices. The conditional inequality is strict, because we found a graph with $8$ vertices which is not linearly free.
It is also proved that for $n, m \leq 6$ the complete bipartite graph $K_{n,m}$ is linearly free.
\end{abstract}

\maketitle



\section{Introduction}
Considering graphs to be topological 1-complexes and thinking over their embeddings into 2 or 3-dimensional spaces, a number of concepts can be established as intrinsic properties of graphs themselves. The minimal genus of graphs can be cited as an example of such properties. A graph is said to be planar (or of genus $0$) if it can be embedded into the Euclidean 2-space $\mathbb{R}^2$. Kuratowski's Theorem gives a complete characterization of planar graphs.
The intrinsical linkedness and the intrinsical knottedness are mentioned as
examples of properties which are derived from embedding into the 3-space.
A graph is said to be intrinsically linked (resp. intrinsically knotted) if every its embedding into the Euclidean 3-space $\mathbb{R}^3$ contains a non-splittable link as a pair of disjoint cycles (resp. a non-trivial knot as a cycle). Otherwise it is said to be linkless (resp. knotless) \cite{CG}.
Linkless graphs were completely characterized by the work of Robertson, Seymour and Thomas \cite{RST}. Their work attracted attention from scholars in low dimensional topology, because a column supporting their proof stands on some firm results from low dimensional topology \cite{ST, Wu, Go}. Currently the characterization of knotless graphs seems to be far from completion, although numerous minimal forbidden graphs have been discovered until recent years \cite{F1, F2, FL, FMN, GMN, HNT, KMO, LKLO, MR}. In this paper we study another intrinsic property which is derived from embeddings of simple graphs.

\begin{figure}[h]
\centering
\includegraphics[width=8cm]{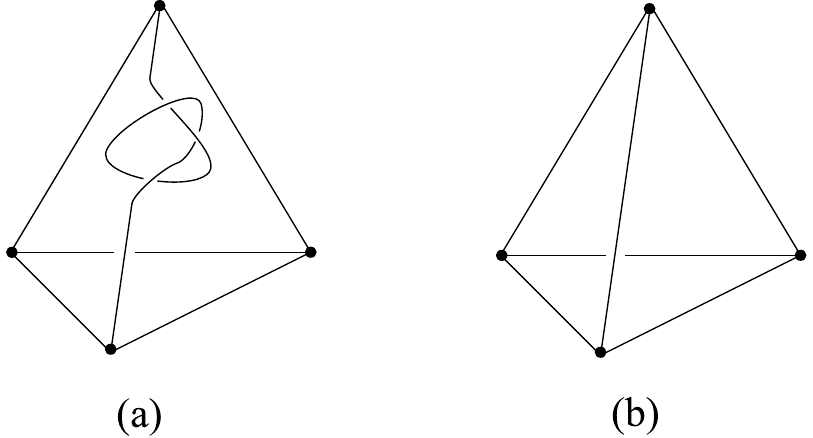}
\caption{Two embeddings $f$ and $g$ of $K_4$}
\label{fig1}
\end{figure}

To address the motivation of our study we observe two specific embeddings $f$ and $g$ of the complete graph $K_4$ into $\mathbb{R}^3$ which are illustrated in Figure \ref{fig1}-(a) and (b), respectively. In Figure \ref{fig1}-(a) an edge forms a locally knotted arc. We see that the fundamental group $\pi_1(\mathbb{R}^3-f(K_4))$ is not free, because it contains the trefoil knot group as a free factor.  For any graph with a cycle, we may construct such a non-free embedding in the same way.
On the contrary if we embed $K_4$ so that every edge is a line segment as seen in Figure \ref{fig1}-(b), then the embedded graph $g(K_4)$ constitutes the 1-skeleton of a tetrahedron, and the fundamental group of its complement is free in consequence. So it can be said that the freeness of such {\em linear} embeddings is an intrinsic property of $K_4$. In fact every complete graph $K_n$ has such property (See Theorem \ref{thmA}).
Motivated by this observation we are interested in discovering more graphs with the property.

\begin{figure}[h]
\centering
\includegraphics[width=8cm]{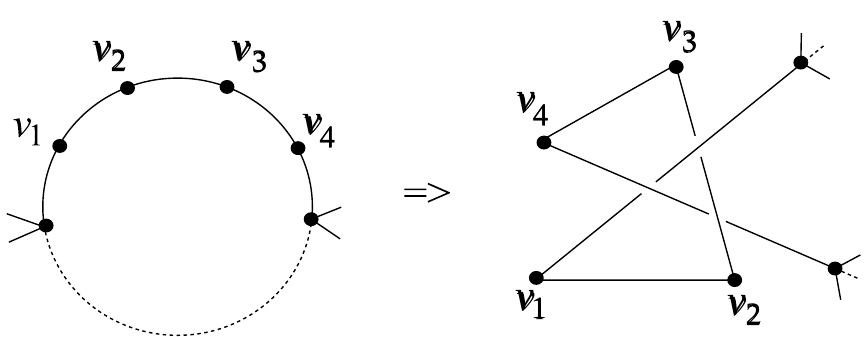}
\caption{A locally knotted arc}
\label{fig2}
\end{figure}

Before stating the main result of this paper we give one more example.
Let $G$ be a simple connected graph. Suppose that a cycle of $G$ contains $n$ consecutive vertices $v_1, \ldots, v_n$ ($n\geq 4$) such that the valency of every $v_i$ is 2. Then we can construct an embedding of $G$ so that every edge is a line segment, the path traversing the $n$ vertices forms a locally knotted arc in a small 3-ball and the 3-ball is disjoint from the other part of the embedded graph. In consequence the fundamental group of the complement is not free (See Figure \ref{fig2} for example).
To avoid some difficulties coming from the flexibility of embedded graph such as local knotting, our study will focus mainly on graphs with valency $\geq 3$ at each vertex.

Throughout this paper graphs are considered to be topological 1-complexes.
\begin{definition} \mbox{}
\begin{enumerate}
\item An embedding of a graph into $\mathbb{R}^3$ is said to be {\em linear}\footnote{ Note that, for a graph to be linearly embeddable, it should be simple, that is, have no multiple edges between any two vertices and no loop edges.} if every embedded edge is a line segment.
\item An embedding $f:H \rightarrow \mathbb{R}^3$ of a graph $H$  is said to be {\em free} if the fundamental group $\pi_1(\mathbb{R}^3-f(H))$ is a free group.
\item A simple connected graph is said to be {\em linearly free} if every linear embedding of the graph is free.
\end{enumerate}
\end{definition}
\noindent \textbf{Notation.} For a graph $G$, let $V(G)$ and $E(G)$ denote its vertex set and edge set, respectively. The degree or valency $deg(v)$ of a vertex $v$ is the number of edges at $v$. And the minimum degree (or valency) of $G$ is denoted by $\delta(G)$, that is, $\delta(G)=\mbox{Min}\;\{deg(v) \; | \; v \in V(G) \}$.
The {\em clique number} $\omega(G)$ of $G$ is the greatest number $r$ such that the complete graph $K_r$ is a subgraph of $G$.

\vspace{0.3cm} In the late 1980s, V. Nicholson proved the following theorem which says
that every complete graph is linearly free.
\begin{ctheorem}{\rm \cite{Nicholson}} \label{thmA}
Every linear embedding of the complete graph $K_n$ is free.
\end{ctheorem}
On the contrary the authors of this paper showed the following theorem:
\begin{ctheorem}{\rm \cite{HL}} \label{thmB}
\mbox{}
\begin{enumerate}
\item  For any $k \geq 1$, there are infinitely many simple connected graphs with minimal valency $k$ which are not linearly free.
\item For any $k \geq 1$, there are infinitely many $k$-connected graphs which are not linearly free.
\end{enumerate}
\end{ctheorem}
\noindent The proof of Theorem \ref{thmB} in \cite{HL} is constructive. The authors constructed two families of graphs which satisfy the two statements respectively. For (1), the constructed graphs have at least $6(k+1)$ vertices, and for (2), at least $12k$ vertices. Note that for the complete graph $K_n$, the number of vertices, the minimal valency and the connectivity index are $n$, $n-1$ and $n-1$, respectively. This indicates that even a graph with high minimal valency or high connectivity may not be linearly free, when the number of vertices is much higher.

\begin{figure}[h]
\centering
\includegraphics[width=5cm]{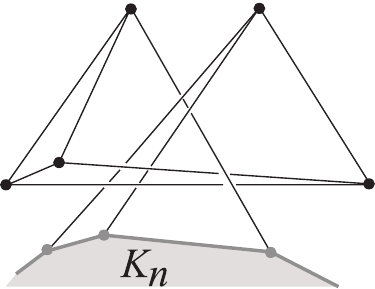}
\caption{$K_n$ with additional five vertices}
\label{fig3}
\end{figure}

\begin{figure}[h]
\centering
\includegraphics[width=8cm]{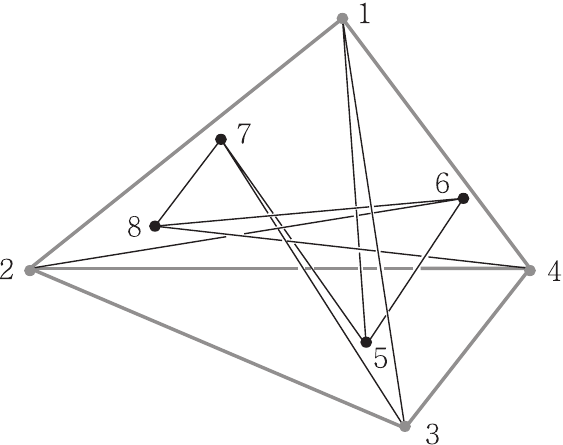}
\caption{$K_4$ with additional four vertices}
\label{fig4}
\end{figure}

\vspace{0.3cm} In this paper we formulate a notion called {\em descending direction} which is a generalization of Nicholson's proof of Theorem \ref{thmA}, and establish a sufficient condition for an embedding of a graph to be free:
\begin{theorem} \label{thm1}
If a linear embedding of a simple graph into $\mathbb{R}^3$ has a descending direction, then it is free.
\end{theorem}
\noindent And applying the condition we find more linearly free graphs. Motivated by Theorem \ref{thmA} we may raise a question: {\em Starting with a complete graph, how much can we enlarge it so that the linear-freeness is not broken and the clique number does not increase?} \\
\noindent The following theorem gives a partial answer.
\begin{theorem} \label{thm2}
Let $G$ be a simple connected graph which contains $K_n$ as a subgraph. If $|V(G)| \leq n+3$, then $G$ is linearly free.
\end{theorem}
\noindent
We remark that the linear-freeness may be broken for graphs with $\delta\geq 3$ and $\omega=n$, when $|V|$ is more than $n+3$. Figure \ref{fig3} depicts a linear embedding of a graph, where the graph is obtained by adding 5 vertices into $K_n$. In the embedding the additional vertices are placed outside the convex hull of the subgraph $K_n$. In consequence the fundamental group of the complement is not free, because it has the trefoil knot group as a free factor. Figure \ref{fig4} illustrates a linear embedding of another graph, where the graph is obtained by adding $4$ vertices into $K_4$. It will be shown that the embedding is not free in a later section.

\begin{figure}[h]
\centering
\includegraphics[width=8cm]{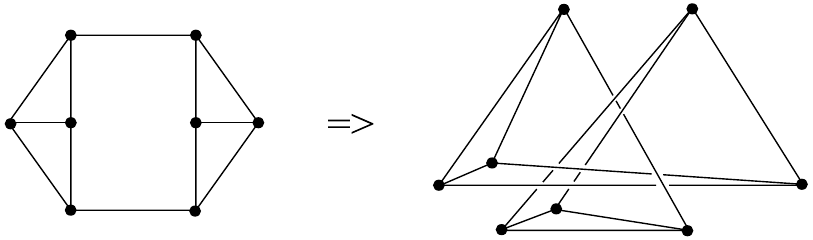}
\caption{A simple connected graph with $|V|=8$ and $\delta=3$}
\label{fig5}
\end{figure}

We are also interested in the linearly free graphs with small number of vertices. A question regarding this is the maximal number of vertices which guarantees the linear-freeness in general for graphs with $\delta \geq 3$. Figure \ref{fig5} illustrates a linear embedding of a graph with $8$ vertices. The fundamental group of the complement is isomorphic to the free product of four infinite cyclic groups and the trefoil knot group. Therefore the following theorem gives the answer for the question.
\begin{theorem} \label{thm3}
 Every simple connected graph with $|V|\leq 7$ and $\delta\geq 3$ is linearly free.
\end{theorem}
The final result is the linear-freeness of the complete bipartite graphs with small number of vertices.
\begin{theorem} \label{thm4}
For $n, m \leq 6$ the complete bipartite graph $K_{n,m}$ is linearly free.
\end{theorem}
\noindent It is not known yet whether every complete bipartite graph is linearly free, which would be one of our possible future works.

\vspace{0.3cm} In Section 2 we introduce the notion of descending direction.
And Theorem \ref{thm1}--\ref{thm4} are proved in the succeeding sections. In the final section it is shown that the embedding in Figure \ref{fig4} is not free.

Throughout the rest of this paper we assume that linear embeddings of graphs are in general position, that is, no three vertices of an embedded graph in $\mathbb{R}^3$ are collinear. Furthermore, when we consider an orthogonal projection from $\mathbb{R}^3$ to a plane, it is also assumed that no three projected vertices are collinear. These assumptions are established by slightly perturbing the embedded graph in $\mathbb{R}^3$.
\section{Descending Direction}
Let $G$ be a simple graph and $f:G\rightarrow \mathbb{R}^n$ be a linear embedding, where $n=2$ or $3$. For our convenience, rather than writing $f(G)$, $f(v)$ and $f(e)$, we allow $G$, $v$ and $e$ to denote the embedded graph, an embedded vertex and an embedded edge, respectively.

For a vertex $v$ of $G$, let $w$ be another vertex which is connected to $v$ by an edge. Then the vector $\vv{vw}$ will be called a {\em $G$-vector} at $v$.
A {\em direction} is a unit vector in $\mathbb{R}^n$. And a direction is said to be {\em generic} with respect to $G$, if it is not orthogonal to any $G$-vector.

For a generic direction $\mathbf{l}$, a vertex $v$ is called a {\em descendant} along $\mathbf{l}$, if there exists a $G$-vector $\vv{vw}$ at $v$ such that
$\mathbf{l}\cdot \vv{vw}<0$. Finally a generic direction $\mathbf{l}$ will be called a {\em descending direction} of the embedded graph $G$, if every vertex of $G$, except for only one, is a descendant along $\mathbf{l}$.
Note that for any generic direction there exists at least one vertex of $G$ which is not a descendant along the direction. One of such vertices belongs to the boundary of the convex hull of $V(G)$.

\section{Proof of Theorem \ref{thm1}}\label{sec-thm1}
Throughout this section  let $G$ be a simple graph and $f:G\rightarrow \mathbb{R}^3$ be a linear embedding with a descending direction. Then we have to prove that $\pi_1(\mathbb{R}^3-f(G))$ is a free group.

By rigid motions of $f(G)$ in $\mathbb{R}^3$ we may assume the followings:
\begin{itemize}
\item The vector $\langle 1,0,0 \rangle$ is a descending direction of $f(G)$.
\item Let $\pi: \mathbb{R}^3 \rightarrow \mathbb{R}^2\times\{0\}$ be the orthogonal projection and $\overline{G}=(\pi \circ f)(G)$. Then the map $\pi \circ f: G \rightarrow \overline{G}$ is {\em regular}, that is, every multiple point of $\overline{G}$ is a transversal double point which is away from vertices.
\item Consider $\mathbb{R}^2\times\{0\}$ to be $\mathbb{R}^2$. Then $\overline{G}$ is contained in the right half-plane $\{(x,y)\in \mathbb{R}^2\;|\; x\geq 0 \}$. Note that the descending direction $\langle 1,0,0 \rangle$ is projected to be the vector $\langle 1, 0 \rangle$ in $\mathbb{R}^2$ by $\pi$.
\end{itemize}

\begin{figure}[h]
\centering
\includegraphics[width=9cm]{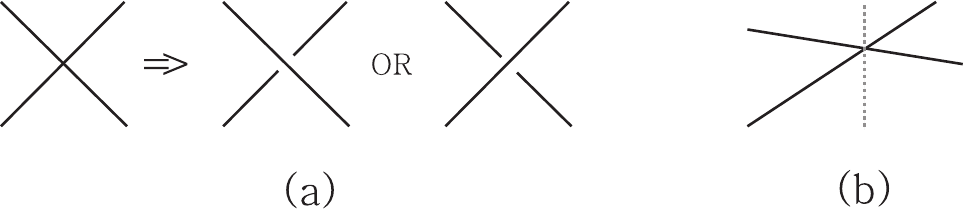}
\caption{(a) Projection and diagram (b) $p$-vertex}
\label{fig6}
\end{figure}

Now, to each double point of $\overline{G}$, add the information on which edge passes over or under in $\mathbb{R}^3$ as illustrated in Figure \ref{fig6}-(a). Then we obtain a {\em diagram} $\tilde{G}$ which represents the embedding $f(G)$.

From now on we consider $\overline{G}$ to be a planar graph on $\mathbb{R}^2$ by taking the vertices of $G$ and the double points of $\overline{G}$ as its vertices. The former will be called {\em $g$-vertices} and the latter {\em $p$-vertices}. Then the projected vector $\mathbf{l}=\langle 1, 0 \rangle$ is a descending direction of the embedded graph $\overline{G}$ on $\mathbb{R}^2$. Recalling that the vector $\langle 1, 0, 0 \rangle$ is a descending direction of $f(G)$ we see that every $g$-vertex of $\overline{G}$, except for the left-most one, is a descendant along $\mathbf{l}$. If $w$ is a $p$-vertex, then it is actually the intersection of two line segments. Since the vector $\mathbf{l}$ is generic, two of the four edges at $w$ are in the left side with respect to a line orthogonal to $\mathbf{l}$ (See Figure \ref{fig6}-(b)). So the vertex $w$ is a descendant of $\overline{G}$ along $\mathbf{l}$.

\vspace{0.3cm}
We will prove Theorem \ref{thm1} in two steps. \\
\noindent {\em Step 1. Construct a spanning tree $T$ of the graph $\overline{G}$}.

\begin{figure}[h]
\centering
\includegraphics[width=7cm]{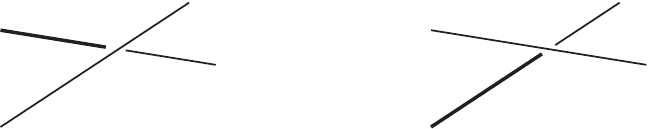}
\caption{Choice of $e_i$ at a $p$-vertex $v_i$}
\label{fig7}
\end{figure}

We may assume that the $x$-coordinates of the vertices of $\overline{G}$ are different with each other. Label the vertices by $v_1$, $v_2$, $\ldots$, $v_{|V(\overline{G})|}$ in the order that their $x$-coordinates increase. Then every $v_i$ except for $v_1$ is a descendant along $\mathbf{l}$, and the left-most two vertices $v_1$ and $v_2$ are $g$-vertices. Now we choose the edges $e_2$, $e_3$, $\ldots$, $e_{|V(\overline{G})|}$ of $T$ as follows. Since the vertex $v_2$ is a descendant there exists an edge, say $e_2$, between $v_2$ and $v_1$. Include $e_2$ into $T$. If $v_i$($i\geq 3$) is a $g$-vertex, then there exists at least one edge which connects $v_i$ to $v_k$ for some $k<i$. Among such edges choose one as $e_i$ so that $k$ is the largest. If $v_i$($i\geq 3$) is a $p$-vertex, then there are two edges which are connected to $v_i$ in its left side. Among the two edges choose one as $e_i$ so that it corresponds to a under-passing at $v_i$ in $\tilde{G}$ (See Figure \ref{fig7}).

\noindent {\em Step 2. Reduce the Wirtinger presentation which is obtained from $\tilde{G}$.}

Let each edge of $\overline{G}$ be oriented so that it goes from the left end vertex to the right. And consider the Wirtinger presentation of $\pi_1(\mathbb{R}^3-f(G))$ which is obtained from $\tilde{G}$ and the orientation of edges. Then each edge of $\overline{G}$ (in fact, a simple loop winding around the edge and based at a fixed point) corresponds to a generator of the presentation. And each vertex $v_i$ of $\overline{G}$ corresponds to a relator, say $W_i$. Denote the generator set by $\mathcal{G}$ and the relator set by $\mathcal{R}=\{W_i\; | \; i=1, 2, \ldots, |V(\overline{G})|\}$. Especially the generator corresponding to the edge $e_i$ of $T$ will be denoted by $x_i$. Now we transform the group presentation by a rewriting process at each $v_i$ from $i=2$ to $|V(\overline{G})|$.

\begin{figure}[h]
\centering
\includegraphics[width=12cm]{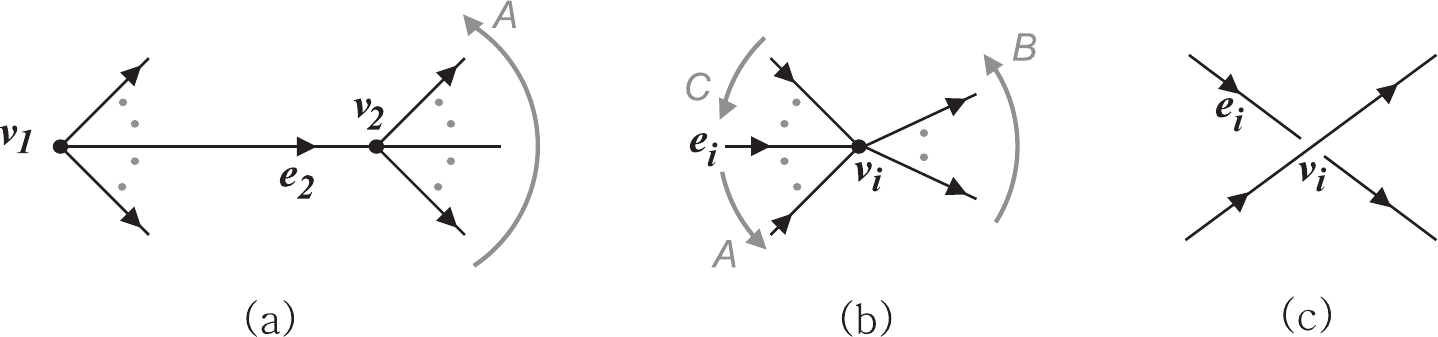}
\caption{Local pictures of $\tilde{G}$ near $v_i$}
\label{fig8}
\end{figure}

Near $v_2$, the local picture of $\tilde{G}$ can be illustrated as in Figure \ref{fig8}-(a). The relator $W_2$ is of the form $W_2=x_2^{-1}A$, where $A$ is a word written in the generators other than $x_2$. Now rewrite $W_1$ by replacing $x_2$ with $A$. Note that $W_1$ and $W_2$ are the only relators in which $x_2$ appears. So we can obtain a new presentation by removing $x_2$ from $\mathcal{G}$, and $W_2$ from $\mathcal{R}$.

In general, if $v_i$ ($i\geq 3$) is a $g$-vertex, then the local picture can be depicted as in Figure \ref{fig8}-(b). The relator $W_i$ would be of the form $W_i=x_i^{-1}ABC$, where $A$, $B$ and $C$ are words written in the generators other than the $i-1$ generators $x_2$, $\ldots$, $x_i$. Figure \ref{fig8}-(c) depicts the local picture in the case that $v_i$ is a $p$-vertex. In this case $W_i$ is of the form $W_i=x_i^{-1}y^{-1}zy$, where $y$ and $z$ are generators other than $x_2$, $\ldots$, $x_i$.

Assume that we have reached the status where the relators $W_2$, $\ldots$, $W_{i-1}$ and the generators $x_2$, $\ldots$, $x_{i-1}$ have been already removed from the presentation. Then $W_i$ and $W_1$ are the only relators in which $x_i$ appears. Now rewrite $W_1$ by substituting $x_i$ with $ABC$ or $y^{-1}zy$. And remove $x_i$ from $\mathcal{G}$, and $W_i$ from $\mathcal{R}$. Repeat this process until $i=|V(\tilde{G})|$. Then, in the final presentation, the generators $x_2$, $\ldots$, $x_{|V(\tilde{G})|}$ are removed and only one relator $W_1$ remains. In fact the resulting $W_1$ is a word obtained by reading the generators that are encountered during going around the boundary curve of a tubular neighborhood of the spanning tree $T$ on $\mathbb{R}^2$. The tubular neighborhood, say $N$, is a homeomorphic $2$-disk, and $\overline{G}-N$ is a disjoint union of simple arcs such that their end-points are on $\partial N$. Therefore we see that $W_1$ consists of some cancellable words of the form $w_k\cdots w_2w_1w_1^{-1}w_2^{-1}\cdots w_k^{-1}$. In conclusion $W_1$ can be removed from the presentation, which implies that $\pi_1(\mathbb{R}^3-f(G))$ is a free group. The proof is completed.

\vspace{0.3cm}
We remark that Theorem \ref{thmA} can be proved as a corollary of Theorem \ref{thm1}. Let $f:K_n \rightarrow \mathbb{R}^3$ be a linear embedding. Then every generic direction is a descending direction of $f(K_n)$, because there exists a straight-line edge between any pair of vertices.

Also Theorem \ref{thm2}, \ref{thm3} and \ref{thm4} will be proved by showing that every linear embedding of the graphs has a descending direction.
\section{Proof of Theorem \ref{thm2}}
We consider the case $|V(G)|=n+3$ in first. Let $H$ be a subgraph of $G$ which is the complete graph $K_n$. And let $u$, $v$ and $w$ be the vertices not belonging to $H$.

Let $f:G \rightarrow \mathbb{R}^3$ be a linear embedding.
As remarked in the end of Section \ref{sec-thm1}, every generic direction is a descending direction of $f(H)$. Therefore it is enough to find a generic direction along which the three vertices are descendants. Without loss of generality we may assume that three linearly independent vectors $\langle 1, 0, 0 \rangle$, $\langle a_1, a_2 ,0 \rangle$ and $\langle b_1, b_2, b_3 \rangle$ are $G$-vectors at $u$, $v$ and $w$, respectively. Then we have to find a vector $\langle x, y, z \rangle$ satisfying
$$ \langle 1,0,0 \rangle \cdot \langle x,y,z \rangle<0, \;\;\;\;
 \langle a_1,a_2,0 \rangle \cdot \langle x,y,z \rangle<0, \;\;\;\;
 \langle b_1,b_2,b_3 \rangle \cdot \langle x,y,z \rangle<0, $$
 that is, $x<0$, $a_1x+a_2y<0$ and $b_1x+b_2y+b_3z<0$. Since this system of inequalities has solutions, there exists a descending direction of $f(G)$.

 When $|V(G)|$ is $n+2$ or $n+1$, we have more freedom to choose a descending direction.

\section{Proof of Theorem \ref{thm4}}
Let $G=K_{n,m}$. Then the vertices of $G$ are split into two families, say {\em white vertices} and {\em black vertices}. Now suppose that $G$ is not linearly free. Then, by Theorem \ref{thm1}, there exists a linear embedding $f:G\rightarrow \mathbb{R}^3$ which has no descending direction. As done in the proof of Theorem \ref{thm1}, consider the regular projection $\overline{G}=\pi \circ f(G)$ to be a planar graph on $\mathbb{R}^2$. Then our assumption implies that $\overline{G}$ has no descending direction on $\mathbb{R}^2$.

Let $\mathbf{l}$ be a generic direction in $\mathbb{R}^2$ with respect to $\overline{G}$.
And we observe the first three $g$-vertices of $\overline{G}$ along the direction.

\vspace{0.3cm}
\noindent {\em Case1: The first two $g$-vertices $v_1$ and $v_2$ have different colors.}

Firstly note that every $p$-vertex is a descendant along $\mathbf{l}$. Since $G$ is a complete bipartite graph, there exists an edge between $v_1$ and $v_2$. Also for any other $g$-vertex $v$ there exists an edge connecting $v$ to $v_1$ or $v_2$. Hence $\mathbf{l}$ is a descending direction of $\overline{G}$.

\vspace{0.3cm}
\noindent {\em Case2: The first two $g$-vertices are of the same color, and the third $g$-vertex $v_3$ has the different color.}

\begin{figure}[h]
\centering
\includegraphics[width=8cm]{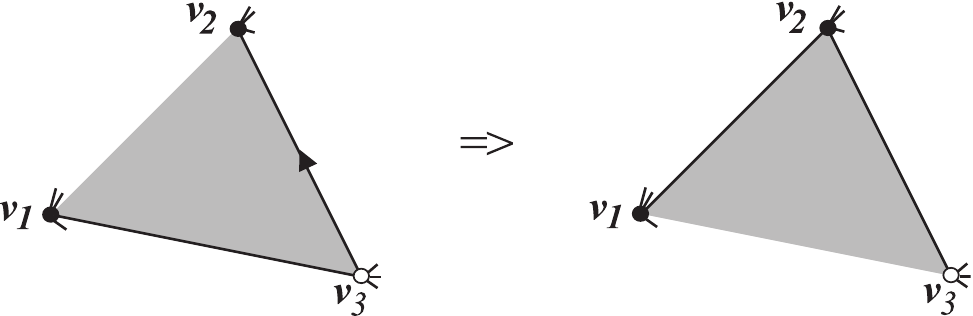}
\caption{Sliding-Move}
\label{fig9}
\end{figure}

In this case any $g$-vertex other than $v_1$ and $v_2$ are descendants along $\mathbf{l}$. Now, for our convenience, let each $v_i$ also denote the corresponding vertex of $f(G)$. Consider the triangle $\Delta_{v_1v_2v_3}$ in $\mathbb{R}^3$ which is determined by the three vertices $v_1$, $v_2$ and $v_3$ of $f(G)$. And modify $f(G)$ by sliding the edge $\overline{v_1v_3}$ along the edge $\overline{v_3v_2}$ as illustrated in Figure \ref{fig9}. We will call this modification a {\em sliding-move}. Let $G^{\prime}$ be the resulting embedded graph in $\mathbb{R}^3$ and $\overline{G^{\prime}}$ be its regular projection.
Then the second $g$-vertex $v_2$ of $\overline{G^{\prime}}$ is connected to $v_1$ by the slid edge. Hence the direction $\mathbf{l}$ is a descending direction of $\overline{G^{\prime}}$, which implies that $\pi_1(\mathbb{R}^3-G^{\prime})$ is free. Note that $f(G)$ has no edge which penetrates the interior of $\Delta_{v_1v_2v_3}$. Therefore the complement $\mathbb{R}^3-f(G)$ is homeomorphic to $\mathbb{R}^3-G^{\prime}$, and
$\pi_1(\mathbb{R}^3-f(G)) \cong \pi_1(\mathbb{R}^3-G^{\prime})$.

\vspace{0.3cm}
By the above observation our proof can be proceeded under an assumption denoted by (\dag) :
{\em For any generic direction, the first three $g$-vertices of $\overline{G}$ have the same color}.

\begin{figure}[h]
\centering
\includegraphics[width=4cm]{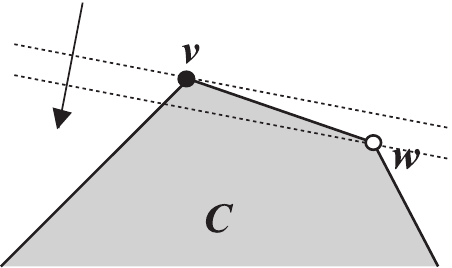}
\caption{Black and white vertices on $P$}
\label{fig10}
\end{figure}

\vspace{0.3cm}
Let $C$ be the convex hull determined by the vertices of $\overline{G}$, and let $P$ be its boundary. Then the vertices of $\overline{G}$ which are not the vertices of the polygon $P$ are included in the interior $\mathring{C}$ of $C$. Note that the vertices of $P$ should have the same color. Otherwise there exist two $g$-vertices $v$ and $w$ of different colors such that the edge $\overline{vw}$ of $\overline{G}$ is also the edge of $P$. Then we can find a generic direction along which $v$ and $w$ are the first two $g$-vertices of $\overline{G}$ as illustrated in Figure \ref{fig10}.

Now, assuming that the vertices of $P$ are black, we observe the necessary number of black vertices in $\mathring{C}$.

\begin{figure}[h]
\centering
\includegraphics[width=11cm]{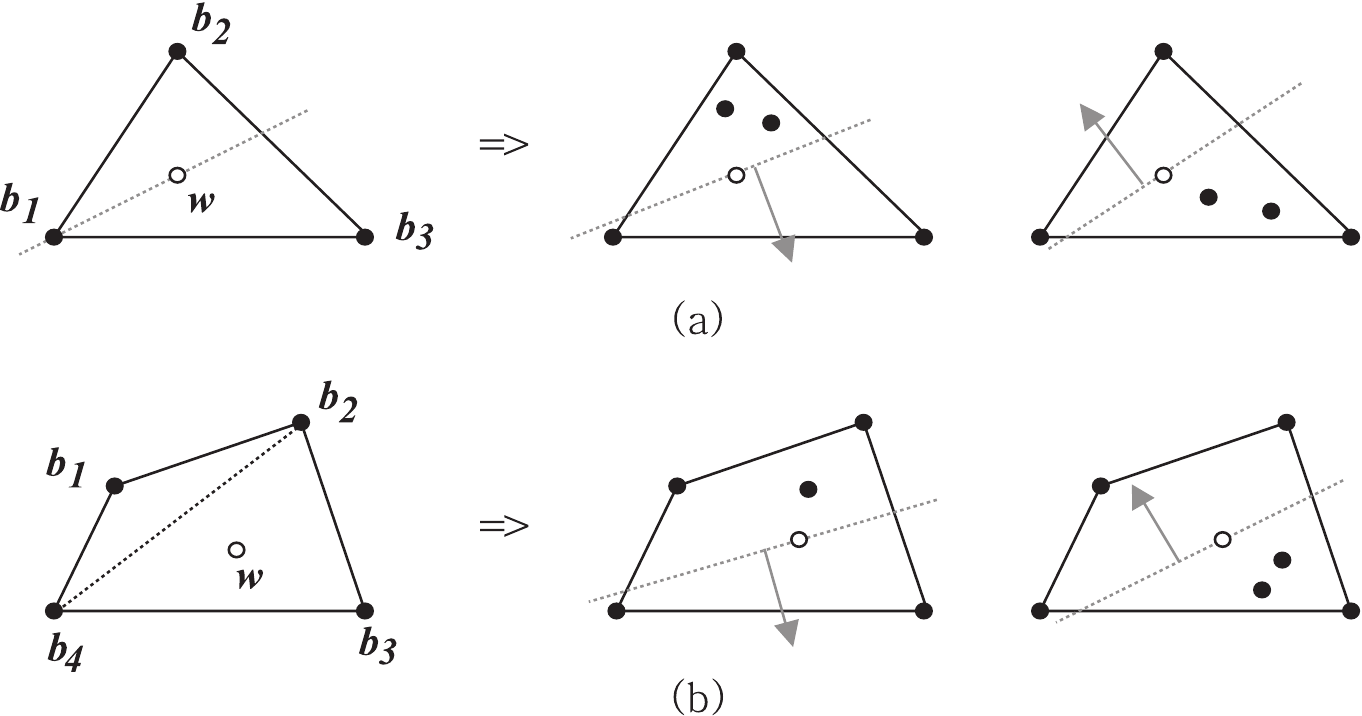}
\caption{Cases: $P$ is a trigon or tetragon.}
\label{fig11}
\end{figure}

(i) Observe the case that $P$ is a trigon with vertices $b_1$, $b_2$ and $b_3$. Select a white vertex $w$. And think over the two directions which are orthogonal to the slight rotations of $\overline{b_1w}$ at $w$ as illustrated in Figure \ref{fig11}-(a). Applying the assumption (\dag) to the two directions we see that $\mathring{C}$ includes at least four black vertices. Therefore $\overline{G}$ has at least seven black vertices.

(ii) In the case that $P$ is a tetragon(with vertices $b_1$, $b_2$, $b_3$ and $b_4$), without loss of generality, we may say that a white vertex $w$ resides inside the trigon $b_2b_3b_4$. See Figure \ref{fig11}-(b). Applying the argument in (i) to $\overline{b_4w}$ we see that $\mathring{C}$ includes at least three black vertices.

\begin{figure}[h]
\centering
\includegraphics[width=11cm]{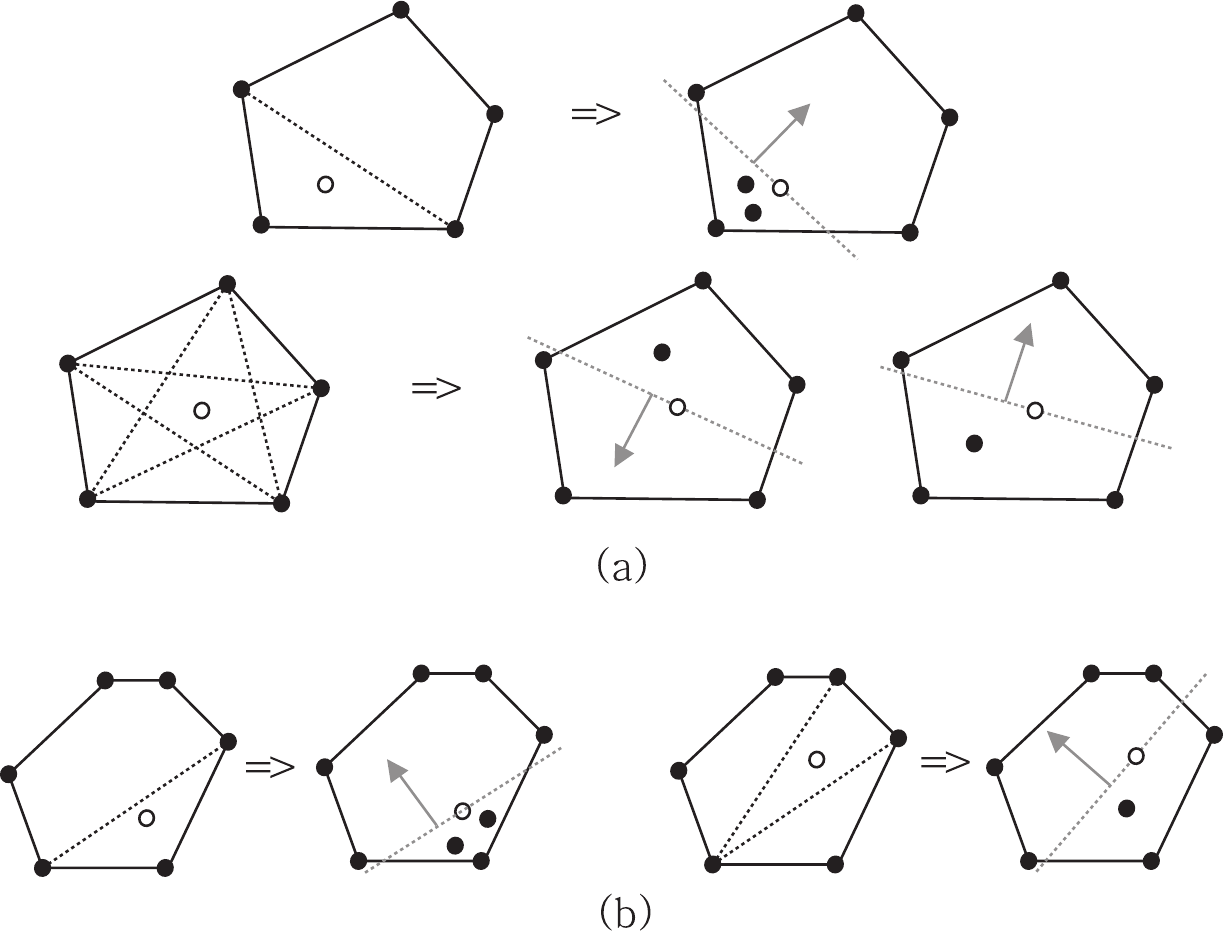}
\caption{Cases: $P$ is a pentagon or hexagon.}
\label{fig12}
\end{figure}

(iii) Similarly Figure \ref{fig12}-(a) and (b) show that if $P$ is a pentagon or a hexagon, then $\mathring{C}$ contains at least two black vertices or one, respectively.

\vspace{0.3cm} By the above observations we know that $G$ has at least seven black vertices. In conclusion, if $K_{n,m}$ has a linear embedding which has no descending direction, then $n\geq 7$ or $m\geq 7$. This deduces Theorem \ref{thm4}.

\section{Proof of Theorem \ref{thm3}}
In this section we will prove Theorem \ref{thm3}. Let $G$ be a simple connected graph with $\delta\geq3$ and $|V(G)|\leq7$.
In a previous work of the authors the following was proved:
\begin{lem}\label{lem1}
\textrm{(Theorem 2 of \cite{HL})}
If $|V(G)|\leq 6$, then $G$ is linearly free.
\end{lem}
Hence, throughout this section, we fix $|V(G)|=7$. For the proof two more lemmas are introduced.
An {$n$-cycle} is a cycle whose number of vertices is $n$.
Lemma \ref{lem2} is easy to prove. Only the proof of Lemma \ref{lem3} is given here.
\begin{lem}\label{lem2}
$G$ contains a $4$-cycle as a subgraph.
\end{lem}
\begin{lem}\label{lem3}
If $G$ contains two disjoint cycles such that one is a $4$-cycle and the other is a $3$-cycle, then every linear embedding of $G$ has a descending direction.
\end{lem}
\begin{proof}
Let $f:G\rightarrow \mathbb{R}^3$ be a linear embedding. For our convenience let $G$ denote the embedded graph $f(G)$ as well as the combinatorial graph $G$ itself. Label the vertices of $G$ by $\{1,2,3,4, 5,6,7\}$, and let $C=(1234)$ be a $4$-cycle of $G$ and $D=(567)$ be a $3$-cycle.
By the condition $\delta \geq 3$ we know that for every vertex $v$ of $D$, there exists a vertex $w$ of $C$ such that $\overline{vw}\in E(G)$. Hence, without loss of generality, it can be assumed that every vertex of $D$ is directly connected to the simple path $\overline{34}\cup \overline{41}$ of $G$ by an edge.

Let $P$ be the plane passing through the three vertices $\{1,3,4\}$. The possible relative positions of the vertices $\{2,5,6,7\}$ with respect to $P$ will be observed. In each case we try to select a descending direction $\mathbf{l}$ between the two unit vectors orthogonal to $P$.

\vspace{0.3cm}
\noindent {\em Case(a): $P$ does not separate $V_1=\{2,5,6,7\}$.}

Select $\mathbf{l}$ so that $\mathbf{l}\cdot \vv{23} < 0$. Then each of $\{2,5,6,7\}$ is a descendant along $\mathbf{l}$. Finally, pushing the vertex $4$ slightly along the counter-direction of $\mathbf{l}$ so that $\mathbf{l}\cdot \vv{34}<0$ and $\mathbf{l}\cdot \vv{14}<0$, the direction $\mathbf{l}$ becomes a descending direction of $G$.

\vspace{0.3cm}
\noindent {\em Case(b): $P$ separates $V_1$ into $\{2\}$ and $\{5,6,7\}$.}

Select $\mathbf{l}$ so that $\mathbf{l}\cdot \vv{23} > 0$. Then each of $\{5,6,7\}$ is a descendant along $\mathbf{l}$. And, pushing the vertex $4$ slightly along $\mathbf{l}$ so that $\mathbf{l}\cdot \vv{43}<0$ and $\mathbf{l}\cdot \vv{41}<0$, the direction $\mathbf{l}$ becomes a descending direction of $G$.

\vspace{0.3cm}
\noindent {\em Case(c): $P$ separates $V_1$ into $\{2,5\}$ and $\{6,7\}$.}

Select $\mathbf{l}$ so that $\mathbf{l}\cdot \vv{23} < 0$.
Then the two vertices $2$ and $5$ are descendants along $\mathbf{l}$. Since there exists an edge between $6$ and $7$, one of the two vertices is a descendant long $\mathbf{l}$. If the vertex $4$ is connected to $6$ or $7$, then push $4$ along the counter-direction of $\mathbf{l}$. Otherwise, one of $1$ and $3$, say $3$, should be connected to $6$ or $7$. In this case push the vertex $4$ slightly and the vertex $1$ a little more along $\mathbf{l}$ so that $\mathbf{l}\cdot \vv{14}<0$ and $\mathbf{l}\cdot \vv{43}<0$ as illustrated in Figure \ref{fig13}.

\vspace{0.3cm}
\noindent {\em Case(d): $P$ separates $V_1$ into $\{2,5,6\}$ and $\{7\}$.}

Select $\mathbf{l}$ so that $\mathbf{l}\cdot \vv{23} < 0$. Then the three vertices $2$, $5$ and $6$ are descendants along $\mathbf{l}$. And, perturbing $\{1,3,4\}$ in the same way with Case(c), the direction $\mathbf{l}$ becomes a descending direction of $G$.
\begin{figure}[h]
\centering
\includegraphics[width=11cm]{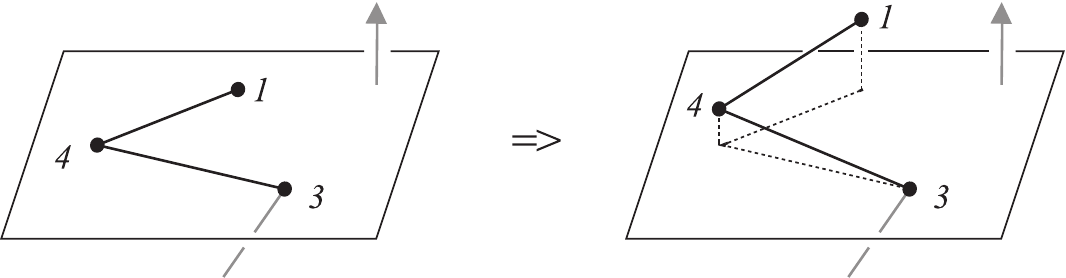}
\caption{Pushing up the vertices $1$ and $4$}
\label{fig13}
\end{figure}
\end{proof}

Now we prove Theorem \ref{thm3}. Let $f:G\rightarrow \mathbb{R}^3$ be a linear embedding.
By Lemma \ref{lem2}, $G$ contains a $4$-cycle. Label the vertices of $G$ by $\{1,2,3,4,5,6,7\}$ so that $(1234)$ is a $4$-cycle. Again this labelling of vertices will be also used to denote the embedded vertices of $f(G)$.

By the condition $\delta\geq 3$, our proof proceeds under an assumption (\ddag): {\em For each vertex $v\in\{5,6,7\}$ there exists an edge which connects $v$ to one of $\{1,3,4\}$}.

Let $P$ be the plane passing through the vertices $\{1,3,4\}$ of $f(G)$.
Now we try to show that a unit vector $\mathbf{l}$ orthogonal to $P$ becomes a descending direction after pushing one or two of $\{1,3,4\}$ as done in the proof of Lemma \ref{lem3}.

If $P$ does not separate $\{5,6,7\}$, then we can find a descending direction of $f(G)$ in the same way with Case(a) or Case(b) of the proof of Lemma \ref{lem3}. If $P$ separates $\{2,5,6,7\}$ into one vertex and three vertices so that the vertex $2$ belongs to the latter, then the argument of Case(d) works. Hence, without loss of generality, it can be assumed that $P$ separates $\{2,5\}$ from $\{6,7\}$. This case is observed in detail.

If there exists an edge between $6$ and $7$, then the argument of Case(c) works, hence the direction $\mathbf{l}$ with $\mathbf{l}\cdot\vv{23}<0$ is selected as a descending direction. Similarly, if there exists an edge between $2$ and $5$, then the direction $\mathbf{l}$ with $\mathbf{l}\cdot\vv{23}>0$ is selected. Consider the case that there is no edge between $6$ and $7$, but there exists a vertex $v\in\{1,3,4\}$ such that both $\overline{v6}$ and $\overline{v7}$ are edges of $f(G)$. In this case the interior of the triangle $\Delta_{67v}$ is disjoint from $f(G)$ in $\mathbb{R}^3$. So we can modify $f(G)$ by sliding $\overline{7v}$ along $v6$, so that $\overline{7v}$ is replaced by $\overline{67}$ in the resulting embedded graph. Then the direction $\mathbf{l}$ with $\mathbf{l}\cdot\vv{23}<0$ is selected.

In conclusion the case that we have to deal with further is described by the following conditions:
\begin{enumerate}
\item[(i)] The plane $P$ separates $\{2,5\}$ from $\{6,7\}$.
\item[(ii)] There is no edge between $2$ and $5$.
\item[(iii)] There is no edge between $6$ and $7$.
\item[(iv)] There is no vertex $v\in\{1,3,4\}$ such that both $\overline{v6}$ and $\overline{v7}$ are edges of $f(G)$.
\end{enumerate}

\vspace{0.3cm}
Applying the assumption (\ddag) and (iv), our case is divided into 12 subcases according to the possible edges of $G$  between $\{1,3,4\}$ and $\{6,7\}$. But, by interchanging $1$ with $3$ and $6$ with $7$, we know that it is enough to investigate the four subcases (S1)$\sim$(S4) in the below. (S1) and (S2) represent the cases that there are exactly three edges between $\{6,7\}$ and $\{1,3,4\}$. (S3) and (S4) represent the cases that there are exactly two edges. For $1 \leq i, j \leq 7$ let $ij$ denote the edge of $G$ between the two vertices $i$ and $j$. And for two subsets $V_1$ and $V_2$ of $V(G)$ let $E(V_1,V_2)$ denote the set of edges of $G$ between $V_1$ and $V_2$.
\begin{enumerate}
\item[(S1)] $E(\{6,7\},\{1,3,4\})=\{61, 64, 73\}$
\item[(S2)] $E(\{6,7\},\{1,3,4\})=\{61, 63, 74\}$
\item[(S3)] $E(\{6,7\},\{1,3,4\})=\{61, 74\}$.
\item[(S4)] $E(\{6,7\},\{1,3,4\})=\{61, 73\}$
\end{enumerate}
In each of the four cases, by the condition $\delta\geq3$, (iii) and (iv), it should be that $E(7, \{2,5,6\})$ is $\{72, 75\}$. Now we look into the subcase (S1). The vertex $6$ should have one more vertex to which it is connected. The possible candidates are $5$ and $2$.

\begin{figure}[h]
\centering
\includegraphics[width=12cm]{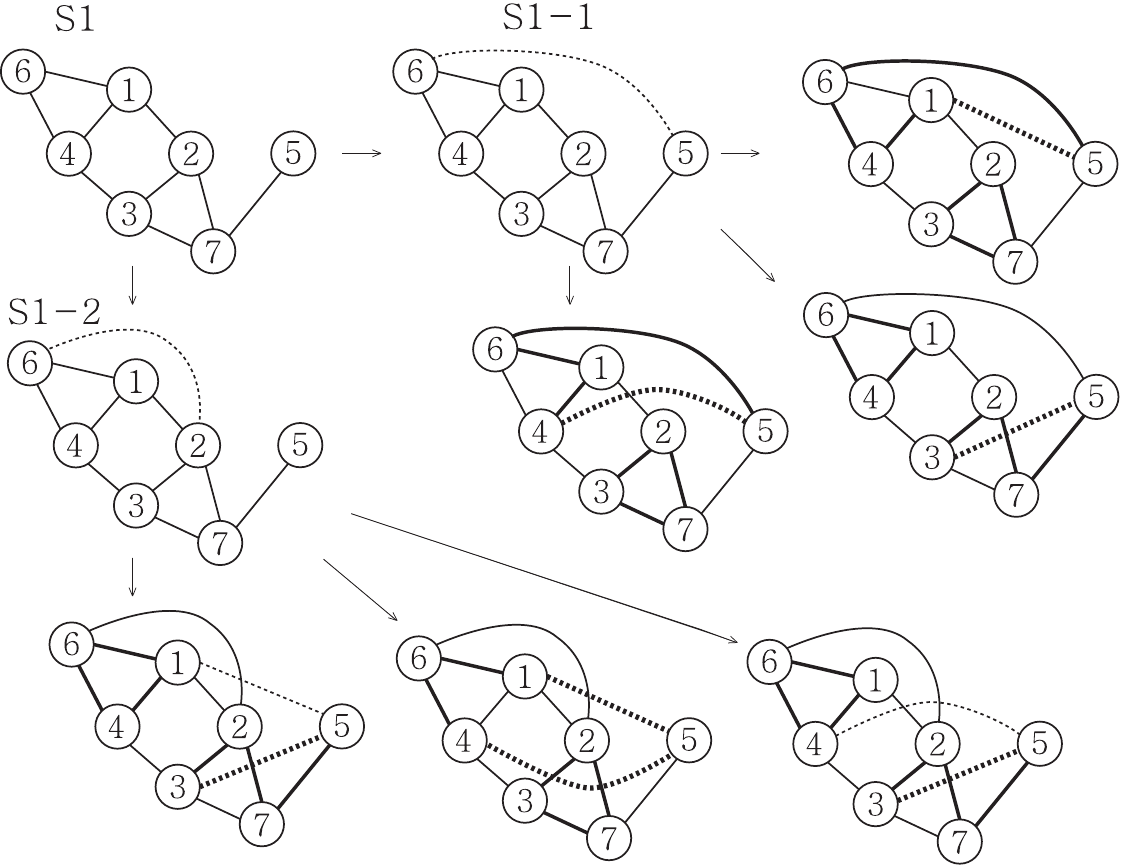}
\caption{Subcase (S1)}
\label{fig14}
\end{figure}

Suppose that the vertex $6$ is connected to $5$ (See (S1-1) in Figure \ref{fig14}). Then also the vertex $5$ should be connected to a vertex $v\in \{1, 3, 4\}$. If $v$ is the vertex $1$, then we can find two disjoint cycles $(1465)$ and $(237)$. If $v$ is $3$ or $4$, then we can find $(2357)$ and $(146)$, or $(1456)$ and $(237)$.

If the vertex $6$ is not connected to $5$, it should be connected to $2$ (See (S1-2) in Figure \ref{fig14}). Then the vertex $5$ needs two more vertices to which it is connected. The possible candidates are $\{1,3\}$, $\{1,4\}$ and $\{3,4\}$. In any of these three possible cases we can find two disjoint cycles $(2357)$ and $(146)$, or $(1546)$ and $(237)$.

\begin{figure}[h]
\centering
\includegraphics[width=12cm]{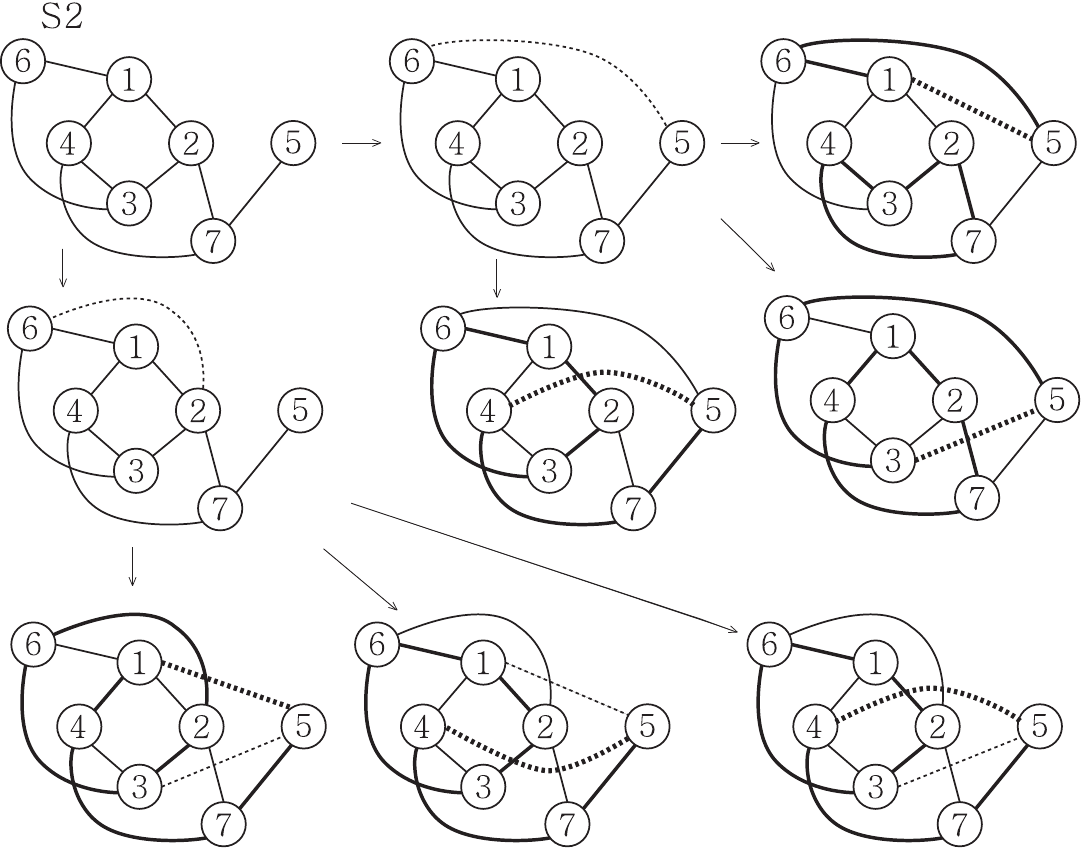}
\caption{Subcase (S2)}
\label{fig15}
\end{figure}

\begin{figure}[h]
\centering
\includegraphics[width=12cm]{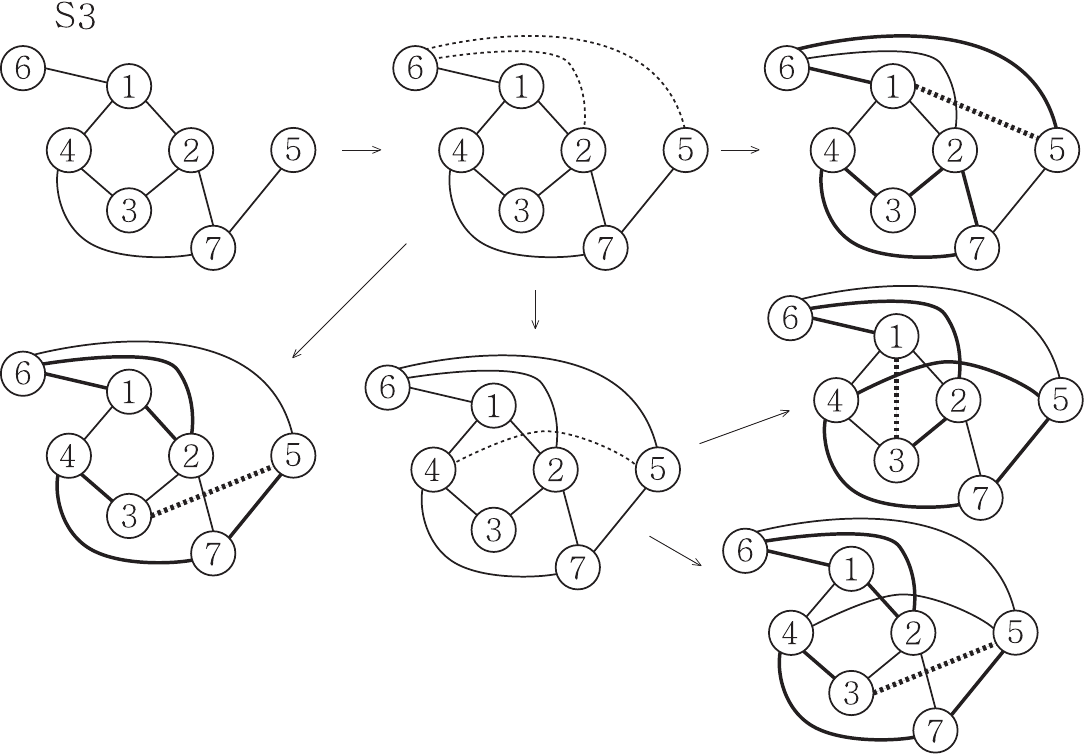}
\caption{Subcase (S3)}
\label{fig16}
\end{figure}

\begin{figure}[h]
\centering
\includegraphics[width=12cm]{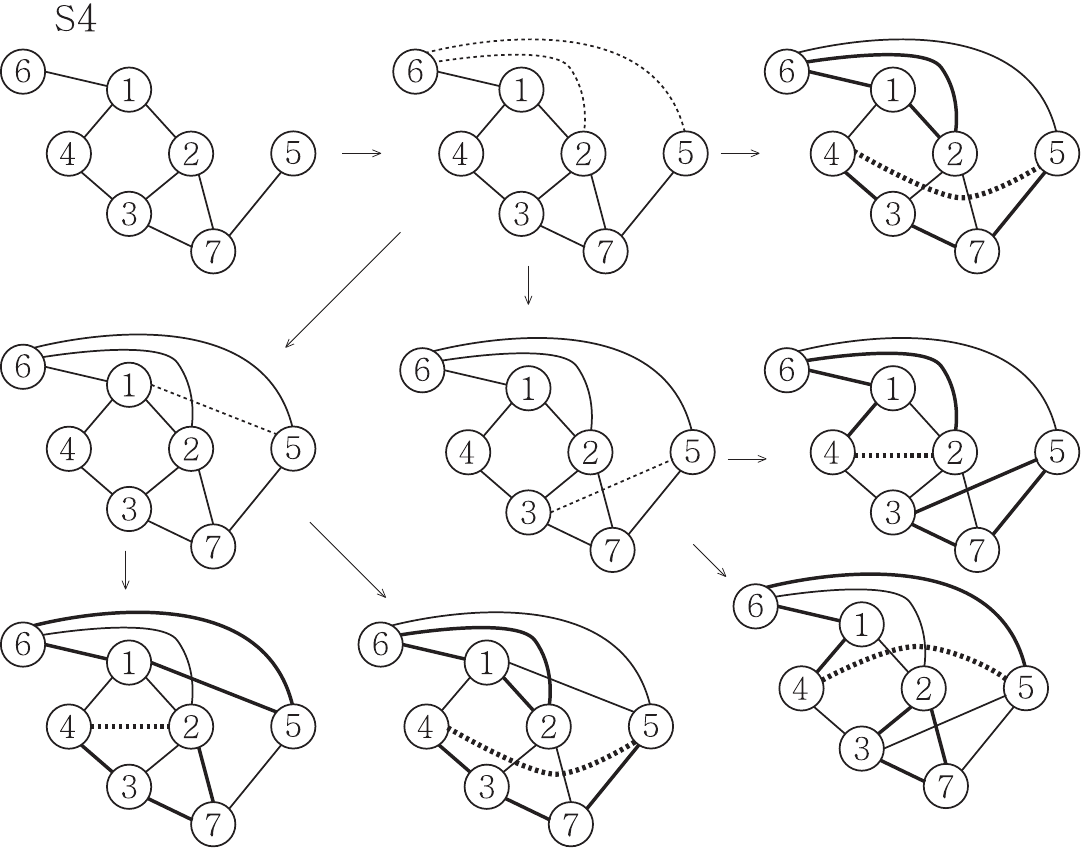}
\caption{Subcase (S4)}
\label{fig17}
\end{figure}

Similarly, for (S2), (S3) and (S4), we can find a $4$-cycle and a $3$-cycle which are disjoint. Figure \ref{fig15}, \ref{fig16} and \ref{fig17} illustrate the edge-recovery process to find the cycles from $G$  in the three cases, respectively. Then, by Lemma \ref{lem3}, the proof of Theorem \ref{thm3} is completed.

\clearpage
\section{Example in Figure \ref{fig4}}
\begin{figure}[h]
\centering
\includegraphics[width=12cm]{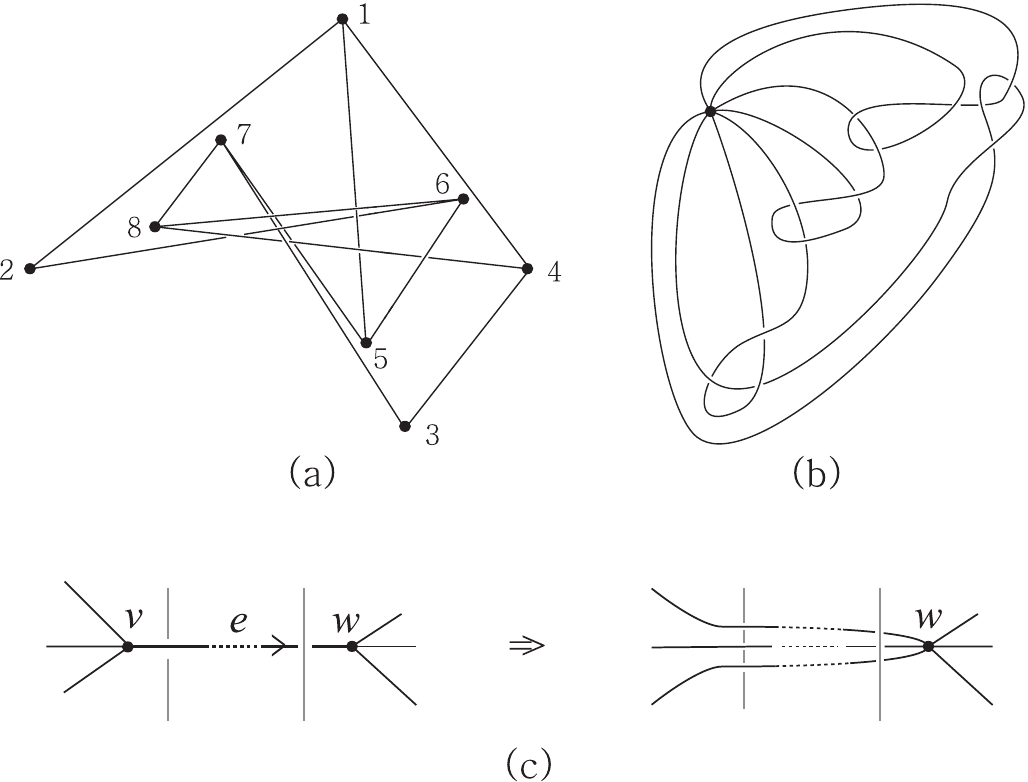}
\caption{$\Gamma'$, $\Theta$, and Edge-contraction $v\rightarrow w$}
\label{fig18}
\end{figure}

Let $\Gamma$ be the embedded graph illustrated in Figure \ref{fig4}.
In this section we show that $\pi_1(\mathbb{R}^3-\Gamma)$ is not a free group.
Let $\Gamma^{\prime}$ be the embedded graph which is obtained from $\Gamma$ by removing three edges $\overline{13}$, $\overline{23}$ and $\overline{24}$ (See Figure \ref{fig18}-(a)). Note that the vertices $5$, $6$, $7$ and $8$ are placed inside the tetrahedron of the other four vertices, hence the interiors of the three triangle $\Delta_{123}$, $\Delta_{124}$ and $\Delta_{134}$ are disjoint from the graph $\Gamma$ itself. This implies that $\pi_1(\mathbb{R}^3-\Gamma^{\prime})$ is a free factor of  $\pi_1(\mathbb{R}^3-\Gamma)$, precisely, $\pi_1(\mathbb{R}^3-\Gamma)\cong \pi_1(\mathbb{R}^3-\Gamma^{\prime})*\mathbb{Z}*\mathbb{Z}*\mathbb{Z}$. So if we show that $\pi_1(\mathbb{R}^3-\Gamma^{\prime})$ is not free, then the proof is completed.

Now we obtain another embedded graph from $\Gamma^{\prime}$ via a sequence of edge- contractions. See Figure \ref{fig18}-(c), where two vertices $v$ and $w$ share an edge $e$. An edge-contraction $v\rightarrow w$ is a contraction of $e$ done by dragging $v$ toward $w$ along the edge. By edge-contractions $1\rightarrow 4$, $4\rightarrow 3$, $3\rightarrow 7$, $8\rightarrow 7$, $5\rightarrow 7$ and $6\rightarrow 7$, and some isotopic moves in $\mathbb{R}^3$, we obtain an embedded bouquet graph $\Theta$ as illustrated in Figure \ref{fig18}-(b).

Let $G$ be a planar graph. An embedding $f:G\rightarrow \mathbb{R}^3$ or the embedded graph $f(G)$ is said to be trivial, if there exists a topological 2-sphere which contains $f(G)$ in $\mathbb{R}^3$.
\begin{lem}\label{lem4}
\textrm{(Theorem 7.5 of \cite{ST})}
Let $G$ be a planar graph. Then an embedding $f:G \rightarrow \mathbb{R}^3$ is trivial if and only if, for every subgraph $H \subseteq G$, $\pi_1(\mathbb{R}^3-f(H))$ is a free group.
\end{lem}

Applying the invariance of 3-coloring for spatial graphs (see Proposition 1.1 in \cite{IY}), we can be aware that the embedded bouquet graph $\Theta$ is not trivial. In fact the diagram of $\Theta$ in Figure \ref{fig18}-(b) is 3-colorable, but the trivially-embedded bouquet graph has no 3-colorable diagram. On the contrary every proper subgraph of $\Theta$ is trivial. Hence Lemma \ref{lem4} implies that $\pi_1(\mathbb{R}^3-\Theta)$ is not free. Note that an edge-contraction preserves the topological type of the complements of the embedded graphs. In conclusion $\pi_1(\mathbb{R}^3-\Gamma^{\prime})$ is not free.

\section*{Acknowledgements}
The first author (corresponding author) was supported by the National Research Foundation of Korea (NRF) grant funded by the Korea government (MSIP) (NRF-2016R1D1A1B01008044).

The second author was supported by the Basic Science Research
Program through the National Research Foundation of Korea (NRF) funded by the Ministry
of Education (2018R1D1A1A09081849).



\end{document}